\documentclass[10pt]{article}

\usepackage{amsmath,amsthm,amssymb,amsfonts}
\usepackage{graphicx} % Required for inserting images
\graphicspath{ {./images/} }
\usepackage{kotex}
\usepackage{authblk} %저자, 소속

\usepackage{enumitem}

\newtheorem{theorem}{Theorem}
\newtheorem*{theorem*}{Theorem}
\newtheorem{lemma}{Lemma}
\newtheorem*{lemma*}{Lemma}

\usepackage[mathlines]{lineno} % peer review, line number

\usepackage{cite}

\usepackage{geometry}

\usepackage{tcolorbox}
\tcbuselibrary{breakable, skins}
\tcbset{enhanced}
\newtcolorbox{tbox}[1]{
    enhanced,
    before skip=2ex, after skip=2ex,
    boxrule=0.5pt, colframe=gray, colback=white, arc=0.5ex,
    boxsep=.5ex, top=1.5ex, bottom=0.5ex, left=0.5em, right=0.5em,
    colbacktitle=white, coltitle=black,
    attach boxed title to top center={xshift=0cm, yshift=-1.5mm},
    boxed title style={size=minimal, enhanced, boxrule=0.25pt, colframe=white},
    fonttitle=\bfseries, breakable=false, title={#1}
}

\title{Range Description for an Attenuated Conical Radon Transform with Fixed Central Axis and Opening Angle}
\author[*]{Gihyeon Jeon}
\affil[*]{School  of  Mathematics,  Kyungpook  National  University,  Daegu  41566,  Republic  of  Korea}
\affil[*]{Corresponding  author: rydbr6709@knu.ac.kr}

\date{}

\begin{document} 
% \linenumbers
\maketitle

\begin{abstract}
    The conical Radon transform is an integral transform that maps a given function $f$ to its integral over a conical surface.  
    In this study, we invesgate the conical Radon transform with a fixed central axis and opening angle, considering the attenuation of radiation within the transform.
    Specifically, we explore the attenuated conical Radon transform.
    In this paper, we provide the range conditions for the attenuated conical Radon transform and its auxiliary transform. 
    Range description of an operator is an important topic in mathematics, and it is useful for understanding the transform, completing incomplete data, improving reconstuction algorithm, correcting measurement error.
    The range conditions of attenuated conical Radon transforms are given in terms of the hyperbolic differential operator.
\end{abstract}
    
Keywords: \texttt{Radon transform, Cone transform, Range description, Compton camera}

\section{Introduction}
The conical Radon transform is an integral transform that assigns to a function $f$ on $\mathbb{R}^{3}$ its surface integral over a cone.
Interest in this transform has increased significantly since the introduction of the Compton camera \cite{singh1983electronically}.
The Compton camera was initially devised for single-photon emission computed tomography. However, it has been applied in various fields, including astronomy, single-scattering optical tomography, and homeland security \cite{terzioglu2018compton}.

Studies on the Compton cameras and various conical Radon transforms are presented in \cite{gouia2014exact, cebeiro2016back, kuchment2016three, moon2019orthogonal, moon2017analytic, nguyen2005radon, truong2007mathematical, kwon2019inversion, jung2016exact, terzioglu2023recovering, moon2020conical, cebeiro2013svd, maxim2009analytical, palamodov2017reconstruction, baines2021range}.
However, most of these studies do not consider radiation attenuation, a crucial factor in image reconstruction.
Neglecting the attenuation effect leads to inaccurate image reconstruction \cite{haltmeier2017inversion}.
We thus consider the attenuated conical Radon transform, which is defined as follows: For $\mu > 0$ and fixed angle $\psi \in (0, \pi/2)$, the attenuated conical Radon transform $C_{\mu}f$ of $f \in C_{c}^{\infty}(\mathbb{R}^{n} \times \mathbb{R})$ is defined by
\begin{align*}
    C_{\mu}f(\mathbf{u}, v) &= \int\limits_{0}^{\infty} \int\limits_{S^{n-1}} f\left( (\mathbf{u}, v) + (r \sin \psi \boldsymbol{\omega}, r \cos \psi)  \right)e^{-\mu r} (r\sin\psi)^{n-1} \mathrm{d}S(\boldsymbol{\omega}) \mathrm{d}r \\ 
    &= \dfrac{\tan^{n-1}\psi}{\cos\psi} \int\limits_{S^{n-1}} \int\limits_{0}^{\infty} f\left( (\mathbf{u}, v) + (z \tan \psi \boldsymbol{\omega}, z) \right)e^{-\frac{\mu}{\cos\psi}z} z^{n-1} \mathrm{d}S(\boldsymbol{\omega}) \mathrm{d}z,
\end{align*}
where $(\mathbf{u}, v) = (u_{1}, \dots, u_{n}, v) \in \mathbb{R}^{n} \times \mathbb{R}$ and the variable is changed by $r\cos \psi \equiv z$.

More research is needed on the attenuated conical Radon transform.
M. Halrmeire et al. proposed the attenuated V-line transform, a two-dimensional variant of the attenuated conical Radon transform, and a reconstruction algorithm using circular detectors \cite{haltmeier2017inversion}.
R. Gouia-Zarrad and S. Moon investigated the inversion formula for the $n$-dimensional attenuated conical Radon transform with a fixed opening angle and provided numerical results \cite{gouia2018inversion}. 
G. Jeon and S. Moon studied the inversion formula of the attenuated conical Radon transform with a fixed opening angle constucted using singular value decomposition \cite{jeon2021singular}.
The generalization of the attenuated conical Radon transform has been further investigated in \cite{haltmeier2018variational, duy2024analysis}.
This study describes the range of the attenuated conical Radon transform and its auxiliary transform.
The range conditions are expressed in terms of the hyperbolic differential operator.

The range of tomographic transforms, such as the Radon transform, is essential in tomographic imaging.
Understanding the range of these transforms is crucial, as it provides significant insight into the measurement data and the underlying functions being reconstructed.
Information about the range of tomographic transforms is useful for several reasons: it aids in comprehending the transforms themselves, completing incomplete data, improving reconstruction algorithms, and correcting measurement errors \cite{agranovsky2007range, kuchment2013radon}.
There have been various studies on the range of different transforms.
G. Ambartsoumian studied the range of the two-dimensional circular Radon transform \cite{ambartsoumian2006range}.
M. Agranovsky et al. provided range conditions for the spherical mean Radon transform \cite{agranovsky2007range}.
A. Katsevich and R. Krylov presented inversion formulas and range conditions for the broken ray transform \cite{katsevich2013broken}.
W. Baines studied the range of conical Radon transform \cite{baines2021range}, and we apply his ideas to the attenuated conical Radon transform.

\begin{figure}[h!]
	\begin{center}
		\begin{tikzpicture}[>=stealth]
		%   \draw (-2,0) -- (4,0) ;
		%   \draw (4,0) -- (6,4) ;
		\draw[->] (0,3) -- (6,3) ;
		\draw[<-] (-2,0) -- (0,3) ;
		\draw[->] (0,3) -- (0,6) ;
		\draw[dashed] (2,1) -- (2,2);
		%\draw[dashed] (2,2) -- (0,4.2);
		\draw[dashed] (2,4.6) -- (3.95,4.6);
		% \draw[dashed] (2,4.6) -- (3.65,4.8);
		\draw[dashed] (-1.3,1) -- (2,1);
		\draw[dashed] (2,1) -- (3.3,3);
		\draw[very thick] (2,2) -- (5,6) ;
		\draw[very thick] (2,2) -- (-1,6) ;
		\draw[->,dashed] (2,2) -- (2,6) ;
		\draw[blue,dashed] (4,4.6) arc (0:180:2cm and 0.3cm);
		\draw[blue] (0,4.6) arc (180:360:2cm and 0.3cm);
		\draw (2,2.5) arc (70:130:10pt);
		\draw[dashed,<->] (2,2) arc (90:270:0.2 and 0.5);
		\draw[<->,dashed,rotate  around={270:(2,2)}] (2,2) arc  (0:180:1.3cm and 0.7cm);
		% \draw (3,4.6) arc (10:30:10pt);
		\node at (2.5,1.8) {$(\mathbf u,v)$};
		\node at (-1.5,1) {$u_1$};
		\node at (3.5,3.2) {$u_2$};
		\node at (3.2,4.6) {$\boldsymbol{\omega}$};
		\node at (1.8,2.7) {$\psi$};
		\node at (6,2.7) {$x_2$};
		\node at (1.7,1.5) {$v$};
		\node at (-1.7,0) {$x_1$};
		\node at (2.8,3.5) {$z$};
		\end{tikzpicture}
	\end{center}
	\caption{A cone $C(\mathbf{u},v)$ of integration for $n=2$ }
	\label{fig:integrationdomain1}
\end{figure}
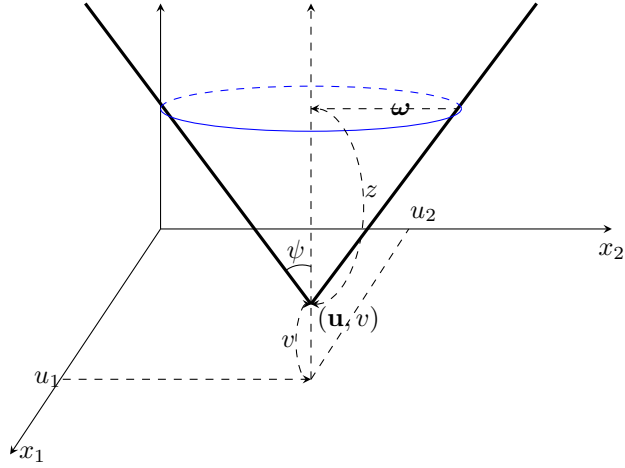

The rest of this paper is structured as follows.
In the next subsection, we present definitions, notations, and lemmas.
Section \ref{sec: main results} outlines the main results regarding the range conditions for the operators $C_{\mu}$. Detailed proofs of these results are provided in Section \ref{sec: proofs}.

\subsection{Preliminary}

The auxiliary operator $A_{\mu}$, analogous to $C_{\mu}$, is useful for characterizing the range of $C_{\mu}$.
The operator $A_{\mu}$ on $C_{c}^{\infty}(\mathbb{R}^{n} \times \mathbb{R})$ is defined similarly as:
\begin{equation*}
    A_{\mu}f(\mathbf{u}, v)
    = \dfrac{\tan^{n-1}\psi}{\cos\psi} \int\limits_{S^{n-1}} \int\limits_{0}^{\infty} f\left( (\mathbf{u}, v) + (z \tan \psi \boldsymbol{\omega}, z) \right)e^{-\frac{\mu}{\cos\psi}z} z^{n-2} \mathrm{d}S(\boldsymbol{\omega}) \mathrm{d}z.
\end{equation*}
In the definition of $A_{\mu}$, the exponent of $z$ is one less than that in $C_{\mu}$.
The Fourier transform of a function $f \in C_{c}( \mathbb{R}^{n})$ is denoted by
\begin{equation*}
    \widehat{f}(\boldsymbol{\xi}) = \mathcal{F}f (\boldsymbol{\xi}) = \int\limits_{\mathbb{R}^{n}} f(\mathbf{x})e^{-\mathrm{i}\boldsymbol{\xi} \cdot \mathbf{x}} \mathrm{d}\mathbf{x}, \qquad \boldsymbol{\xi} = (\xi_{1}, \dots, \xi_{n}).
\end{equation*}
Note that a bounded function $g$ can be considered a tempered distribution, and its Fourier transform exists.
In addition, if smooth function $g$ can be considered as a tempered distribution, then so is $g'$, and the Fourier transform of $g'$ exists because for Schwartz function $\phi$:
\begin{equation*}
    g'[\phi] = - g[\phi'].
\end{equation*}

Linear differential operator $L$, defined as follows, is plays a significant role in describing the range of $C_{\mu}$ and $A_{\mu}$:
\begin{equation*}
    L := \frac{\mu^{2}}{\cos^{2}\psi} - \frac{2\mu}{\cos\psi}\partial_{z} + \partial_{z}^{2} - \tan^{2}\psi \Delta_{\mathbf{x}}.
\end{equation*}
Here, $\Delta_{\mathbf{x}}$ is Laplacian with respect to $\mathbf{x}$.
Since $\widehat{\partial_{j} f}(\boldsymbol{\xi}) = \mathrm{i}\xi_{j}\widehat{f}(\boldsymbol{\xi})$, we have
\begin{equation*}
    \widehat{Lf} (\boldsymbol{\xi}, \sigma) = \left( (\frac{\mu}{\cos\psi} - \mathrm{i}\sigma)^{2} + (\left| \boldsymbol{\xi} \right| \tan\psi)^{2} \right) \widehat{f} (\boldsymbol{\xi}, \sigma),
\end{equation*}
where $(\boldsymbol{\xi}, \sigma)$ is the frequency domain variable corresponding to $(\mathbf{x}, z)$.

We now introduce Lemma needed to prove main Theorems.
The following Lemma is about the Fourier transform of two operators, $C_{\mu}$ and $A_{\mu}$.

\begin{lemma}\label{lem: fourier transform of Cmuf and Amuf}
    For $\mu > 0$ and $(\boldsymbol{\xi}, \sigma) \in \mathbb{R}^{n} \times \mathbb{R}$, we have
    $$
    \widehat{C_{\mu}f} (\boldsymbol{\xi}, \sigma)
    =  \dfrac{ \alpha_{n} \left( \frac{\mu}{\cos\psi} - \mathrm{i}\sigma \right) \widehat{f} (\boldsymbol{\xi}, \sigma)}{\left((\frac{\mu}{\cos\psi} - \mathrm{i}\sigma)^{2} + (\left| \boldsymbol{\xi} \right| \tan\psi)^{2} \right)^{\frac{n+1}{2}}} 
    \quad \text{ and } \quad
    \widehat{A_{\mu}f} (\boldsymbol{\xi}, \sigma) = \dfrac{\beta_{n} \widehat{f} (\boldsymbol{\xi}, \sigma)}{\left( (\frac{\mu}{\cos\psi} - \mathrm{i} \sigma)^{2} + (\left| \boldsymbol{\xi} \right| \tan\psi)^{2} \right)^{\frac{n-1}{2}}} 
    $$
    where $\alpha_{n} = \dfrac{2^{n}\pi^{\frac{n-1}{2}} \Gamma(\frac{n+1}{2}) \tan^{n-1}\psi }{\cos\psi}$, $\beta_{n} = \dfrac{2^{n-1}\pi^{\frac{n-1}{2}}\Gamma(\frac{n-1}{2})\tan^{n-1}\psi}{\cos\psi}$, and $\Gamma$ is the gamma function. 
    
\end{lemma}
Proof of Lemma is given in Appendix. 

\section{Main Results}\label{sec: main results}
We now present the main results of this study.
A necessary and sufficient condition for $g \in C^{\infty}(\mathbb{R}^{n+1})$ is that $L^{k}g$ $(k \in \mathbb{N})$ has compact support.
The range conditions for the operators $C_{\mu}$ are given in the following Theorems.
\begin{theorem}\label{thm: odd dim of Cmu}
    Let $n = 2k-1$ be odd (i.e., $\mathbf{x} \in \mathbb{R}^{2k-1}$). A function $g \in C^{\infty}(\mathbb{R}^{2k-1} \times \mathbb{R})$ is in the range of $C_{\mu}$ on $C_{c}^{\infty}(\mathbb{R}^{2k-1} \times \mathbb{R})$ if and only if the following conditions are satisfied:
    \begin{enumerate}[label=(\roman*)]
        \item $L^{k}g(x, z)$ has compact support.
        \item $\displaystyle \int\limits_{-\infty}^{\infty} e^{-\frac{\mu}{\cos\psi}z}L^{k}g(x, z) \mathrm{d}z = 0$ for every $x \in \mathbb{R}^{2k-1}$.
    \end{enumerate}
\end{theorem}
In the specific case where $n$ is even, the operator $A_{\mu}$ is used to demonstrate the range conditions of $C_{\mu}$.    
\begin{theorem}\label{thm: even dim of Cmu}
    Let $n = 2k$ be even (i.e., $x \in \mathbb{R}^{2k}$). A function $g \in C^{\infty}(\mathbb{R}^{2k} \times \mathbb{R})$ is in the range of $C_{\mu}$ on $C_{c}^{\infty}(\mathbb{R}^{2k} \times \mathbb{R})$ if and only if the following conditions are satisfied:
    \begin{enumerate}[label=(\roman*)]
        \item $L^{2k}A_{\mu}g (x, z)$ has compact support.
        \item $\displaystyle \int\limits_{-\infty}^{\infty} e^{-\frac{\mu}{\cos\psi}z} L^{2k}A_{\mu}g (x, z) \mathrm{d}z = 0$ for every $x \in \mathbb{R}^{2k}$
    \end{enumerate}
\end{theorem}
Theorems regarding range conditions for $A_{\mu}$ are analogous to those for $C_{\mu}$ but are more straightforward, as shown in the following.
\begin{theorem} \label{thm: odd dim of Amu}
    Let $n = 2k-1$ be odd (i.e., $\mathbf{x} \in \mathbb{R}^{2k-1}$). For a function $g \in C^{\infty}(\mathbb{R}^{2k-1} \times \mathbb{R})$, $g$ is in the range of $A_{\mu}$ on $C_{c}^{\infty}(\mathbb{R}^{2k-1} \times \mathbb{R})$ if and only if $L^{k-1}g$ has compact support.
\end{theorem}
In the specific case where $n$ is odd, the range conditions of $A_{\mu}$ are related to $L^{2k-1}A_{\mu}g$, similar to Theorem \ref{thm: even dim of Cmu}.
\begin{theorem} \label{thm: even dim of Amu}
    Let $n = 2k$ be odd (i.e., $\mathbf{x} \in \mathbb{R}^{2k}$). For a function $g \in C^{\infty}(\mathbb{R}^{2k} \times \mathbb{R})$, $g$ is in the range of $A_{\mu}$ on $C_{c}^{\infty}(\mathbb{R}^{2k} \times \mathbb{R})$ if and only if $L^{2k-1}A_{\mu}g$ has compact support.
\end{theorem}

Proofs are provided in the next section.

\section{Proofs of our results} \label{sec: proofs}
This section provides detailed proofs of the Theorems presented in Section \ref{sec: main results}.

\subsection{Proof of Theorem \ref{thm: odd dim of Cmu}}
For $f \in C_{c}^{\infty}(\mathbb{R}^{2k-1} \times \mathbb{R})$, suppose $g = C_{\mu}f$.
By definition of $C_{\mu}$, $g$ is smooth. Since $f$ has compact support, $g$ is bounded and the Fourier transform of $L^{k}g$ is given by
\begin{align*}
    \widehat{L^{k}g} (\boldsymbol{\xi}, \sigma)
    &= \left( \left(\frac{\mu}{\cos\psi} - \mathrm{i} \sigma\right)^{2} + (\left| \boldsymbol{\xi} \right| \tan\psi)^{2} \right)^{k} \widehat{g} (\boldsymbol{\xi}, \sigma) \\
    &= \left( \left(\frac{\mu}{\cos\psi} - \mathrm{i} \sigma\right)^{2} + (\left| \boldsymbol{\xi} \right| \tan\psi)^{2} \right)^{k} \widehat{C_{\mu}f} (\boldsymbol{\xi}, \sigma) \\
    &= \left( \left(\frac{\mu}{\cos\psi} - \mathrm{i} \sigma\right)^{2} + (\left| \boldsymbol{\xi} \right| \tan\psi)^{2} \right)^{k} \alpha_{2k-1} \dfrac{\frac{\mu}{\cos\psi} - \mathrm{i}\sigma}{\left((\frac{\mu}{\cos\psi} - \mathrm{i}\sigma)^{2} + (\left| \boldsymbol{\xi} \right| \tan\psi)^{2} \right)^{k}} \widehat{f} (\boldsymbol{\xi}, \sigma) \\
    &= \alpha_{2k-1}\left( \frac{\mu}{\cos\psi} - \mathrm{i}\sigma \right)\widehat{f} (\boldsymbol{\xi}, \sigma) \\
    &= \alpha_{2k-1}\mathcal{F}\left[ \frac{\mu}{\cos\psi}f - \partial_{z}f \right](\boldsymbol{\xi}, \sigma).
\end{align*}
yielding
\begin{equation*}
    L^{k}g (\mathbf{x}, z) = \alpha_{2k-1} \left( \frac{\mu}{\cos\psi}f - \partial_{z}f \right) (\mathbf{x}, z),
\end{equation*}
or equivalently
\begin{align*}
-e^{-\frac{\mu}{\cos\psi}z} L^{k}g (\mathbf{x}, z) &= \alpha_{2k-1} \left( \frac{-\mu}{\cos\psi} e^{-\frac{\mu}{\cos\psi}z} f + e^{-\frac{\mu}{\cos\psi}z} \partial_{z}f \right) (\mathbf{x}, z)\\
                                    &= \alpha_{2k-1} \dfrac{\partial }{\partial z}(e^{-\frac{\mu}{\cos\psi}z} f) (\mathbf{x}, z).
\end{align*}
Since $f$ has compact support, so is $L^{k}g(\mathbf{x}, z)$.
Besides, we have:
\begin{align*}
    \int \limits_{-\infty}^{\infty}e^{-\frac{\mu \tau}{\cos\psi}} L^{k}g(\mathbf{x}, \tau) \mathrm{d}\tau
    &= \lim\limits_{z \to \infty} \int\limits_{-z}^{z}e^{-\frac{\mu \tau}{\cos\psi}} L^{k}g(\mathbf{x}, \tau) \mathrm{d}\tau \\ 
    &= -\alpha_{2k-1} \lim\limits_{z \to \infty} (e^{-\frac{\mu}{\cos\psi}z} f(\mathbf{x}, z) - e^{\frac{\mu}{\cos\psi}z} f(\mathbf{x}, -z) ) \\ 
    &= 0.
\end{align*}

Conversely, suppose $g \in C^{\infty}$ satisfies conditions {\it(i)} and {\it(ii)}, $f$ can be defined as:
\begin{equation*}
    f(\mathbf{x}, z) := -\alpha_{2k-1}^{-1}e^{-\frac{\mu}{\cos\psi}z} \int\limits_{-\infty}^{z} e^{-\frac{\mu \tau}{\cos\psi}} L^{k}g(x,\tau) \mathrm{d}\tau.
\end{equation*}
Since $g$ is smooth, so is $f$.
Additionally, $f$ has compact support because $g$ satisfies conditions {\it(i)} and {\it(ii)}.
Multiplying both sides by $-e^{\frac{\mu}{\cos \psi}z}$, we obtain:
\begin{equation*}
    -e^{-\frac{\mu}{\cos\psi}z} f(\mathbf{x}, z) = \alpha_{2k-1}^{-1} \int\limits_{-\infty}^{z} e^{-\frac{\mu \tau}{\cos\psi}} L^{k}g(x,\tau) \mathrm{d}\tau.
\end{equation*}
Differentiating this expression with respect to $z$ yields:
\begin{equation*}
    \frac{\mu}{\cos\psi} f(\mathbf{x}, z) - \partial_{z}f(\mathbf{x}, z) = \alpha_{2k-1}^{-1} L^{k}g(\mathbf{x}, z).
\end{equation*}
Taking the Fourier transform, we get:
\begin{equation*}
    \left( \frac{\mu}{\cos\psi} - \mathrm{i}\sigma \right)\widehat{f}(\boldsymbol{\xi}, \sigma) = \alpha_{2k-1}^{-1} \widehat{L^{k}g}(\boldsymbol{\xi}, \sigma).
\end{equation*}
Therefore, the Fourier transform of $C_{\mu}f$ is as follows:
\begin{align*}
    \widehat{C_{\mu}f} (\boldsymbol{\xi}, \sigma) &= \alpha_{2k-1} \dfrac{\frac{\mu}{\cos\psi} - \mathrm{i}\sigma}{\left( (\frac{\mu}{\cos\psi} - \mathrm{i}\sigma)^{2} + (\left| \boldsymbol{\xi} \right| \tan\psi)^{2} \right)^{k}}\widehat{f} (\boldsymbol{\xi}, \sigma)\\
    &= \alpha_{2k-1} \dfrac{\alpha_{2k-1}^{-1}}{\left( (\frac{\mu}{\cos\psi} - \mathrm{i}\sigma)^{2} + (\left| \boldsymbol{\xi} \right| \tan\psi)^{2} \right)^{k}}\widehat{L^{k}g} (\boldsymbol{\xi}, \sigma)\\
    &= \widehat{g} (\boldsymbol{\xi}, \sigma).
\end{align*}

\subsection{Proof of Theorem \ref{thm: even dim of Cmu}}
Proofs of Theorems \ref{thm: even dim of Cmu} and \ref{thm: odd dim of Cmu} are similar.
For $f \in C_{c}^{\infty}(\mathbb{R}^{2k} \times \mathbb{R})$, suppose $g = C_{\mu}f$.
Since $f$ has compact support, $g$ is bounded, and the Fourier transform of $L^{2k}A_{\mu}g$ is given by:
\begin{align*}
    \mathcal{F}[L^{2k}A_{\mu}g] (\boldsymbol{\xi}, \sigma)
    &= \left( \left(\frac{\mu}{\cos\psi} - \mathrm{i} \sigma\right)^{2} + (\left| \boldsymbol{\xi} \right| \tan\psi)^{2} \right)^{2k} \widehat{A_{\mu}g} (\boldsymbol{\xi}, \sigma) \\
    &= \left( \left(\frac{\mu}{\cos\psi} - \mathrm{i} \sigma\right)^{2} + (\left| \boldsymbol{\xi} \right| \tan\psi)^{2} \right)^{2k} \dfrac{\beta_{2k}}{\left( \left(\frac{\mu}{\cos\psi} - \mathrm{i} \sigma\right)^{2} + (\left| \boldsymbol{\xi} \right| \tan\psi)^{2} \right)^{\frac{2k-1}{2}}} \widehat{g} (\boldsymbol{\xi}, \sigma) \\
    &= \left( \left(\frac{\mu}{\cos\psi} - \mathrm{i} \sigma\right)^{2} + (\left| \boldsymbol{\xi} \right| \tan\psi)^{2} \right)^{\frac{2k+1}{2}} \dfrac{\alpha_{2k}\beta_{2k}(\frac{\mu}{\cos\psi} - \mathrm{i}\sigma)}{\left( \left(\frac{\mu}{\cos\psi} - \mathrm{i} \sigma\right)^{2} + (\left| \boldsymbol{\xi} \right| \tan\psi)^{2} \right)^{\frac{2k+1}{2}}} \widehat{f}(\boldsymbol{\xi}, \sigma) \\
    &= \alpha_{2k}\beta_{2k} \left(\frac{\mu}{\cos\psi} - \mathrm{i}\sigma\right) \widehat{f} (\boldsymbol{\xi}, \sigma) \\
    &= \alpha_{2k}\beta_{2k}\mathcal{F} \left[ \frac{\mu}{\cos\psi}f - \partial_{z}f \right] (\boldsymbol{\xi}, \sigma).
\end{align*}
yielding
\begin{equation*}
    L^{2k}A_{\mu}g (\mathbf{x}, z) = \alpha_{2k}\beta_{2k} \left( \frac{\mu}{\cos\psi} f - \partial_{z}f \right) (\mathbf{x}, z),
\end{equation*}
or equivalently
\begin{align*}
    -e^{-\frac{\mu}{\cos\psi}z} L^{2k}A_{\mu}g (\mathbf{x}, z) &= \alpha_{2k}\beta_{2k} \left( -\frac{\mu}{\cos\psi}e^{-\frac{\mu}{\cos\psi}z} f + e^{-\frac{\mu}{\cos\psi}z} \partial_{z}f \right) (\mathbf{x}, z) \\
                                                               &= \alpha_{2k}\beta_{2k} \dfrac{\partial }{\partial z}(e^{-\frac{\mu}{\cos\psi}z} f) (\mathbf{x}, z).
\end{align*}
Since $f$ has compact support, so is $L^{k}A_{\mu} g$. Besides, we find:
\begin{align*}
    \int\limits_{-\infty}^{\infty} e^{-\frac{\mu}{\cos\psi}z} L^{2k}A_{\mu}g (\mathbf{x},\tau) \mathrm{d}\tau
    &= \lim\limits_{z \to \infty} \int\limits_{-z}^{z} e^{-\frac{\mu}{\cos\psi}z} L^{2k}A_{\mu}g (\mathbf{x},\tau) \mathrm{d}\tau \\
    &= -\alpha_{2k}\beta_{2k}\lim\limits_{z \to \infty}(e^{-\frac{\mu}{\cos\psi}z} f(\mathbf{x}, z) - e^{\frac{\mu}{\cos\psi}z} f(\mathbf{x}, -z) ) \\
    &= 0.
\end{align*}

Conversely, suppose $g \in C^{\infty}$ satisfies conditions {\it (i)} and {\it (ii)}, $f$ can be defined as:
\begin{equation*}
    f(\mathbf{x}, z) := -\alpha_{2k}^{-1}\beta_{2k}^{-1}e^{\frac{\mu}{\cos\psi}z} \int\limits_{-\infty}^{z} e^{-\frac{\mu}{\cos\psi}z} L^{2k}A_{\mu}g (\mathbf{x},\tau) \mathrm{d}\tau.
\end{equation*}
Since $g$ is smooth, so is $f$.
Additionally, $f$ has compact support because $g$ satisfies conditions {\it (i)} and {\it (ii)}.
Multiplying both sides by $-e^{-\frac{\mu}{\cos\psi}z}$, we get:
\begin{equation*}
    -e^{-\frac{\mu}{\cos\psi}z}f(\mathbf{x}, z) = \alpha_{2k}^{-1}\beta_{2k}^{-1} \int\limits_{-\infty}^{z} e^{-\frac{\mu \tau}{\cos\psi}} L^{2k}A_{\mu}g (\mathbf{x},\tau) \mathrm{d}\tau,
\end{equation*}
and differentiating with respect to $z$, we have:
\begin{equation*}
    \frac{\mu}{\cos\psi} f(\mathbf{x}, z) - \partial_{z}f(\mathbf{x}, z) = \alpha_{2k}^{-1}\beta_{2k}^{-1} L^{2k}A_{\mu}g (\mathbf{x}, z).
\end{equation*}
Taking the Fourier transform with respect to $(\mathbf{x}, z)$ yields:
\begin{equation*}
   \left( \frac{\mu}{\cos\psi} - \mathrm{i}\sigma \right)\widehat{f}(\xi, \sigma) = \alpha_{2k}^{-1}\beta_{2k}^{-1} \mathcal{F}[L^{2k}A_{\mu}g](\xi, \sigma).
\end{equation*}
Therefore, the Fourier transform of $C_{\mu}f$ is as follows:
\begin{align*}
    \widehat{C_{\mu}f} (\boldsymbol{\xi}, \sigma)
    &= \alpha_{2k} \dfrac{\frac{\mu}{\cos\psi} - \mathrm{i}\sigma}{\left( (\frac{\mu}{\cos\psi} - \mathrm{i}\sigma)^{2} + (\left| \boldsymbol{\xi} \right| \tan\psi)^{2} \right)^{\frac{2k+1}{2}}}\widehat{f} (\boldsymbol{\xi}, \sigma) \\
    &= \alpha_{2k} \dfrac{\alpha_{2k}^{-1}\beta_{2k}^{-1}}{\left( (\frac{\mu}{\cos\psi} - \mathrm{i}\sigma)^{2} + (\left| \boldsymbol{\xi} \right| \tan\psi)^{2} \right)^{\frac{2k+1}{2}}}\widehat{L^{2k}A_{\mu}g} (\boldsymbol{\xi}, \sigma) \\
    &= \dfrac{\beta_{2k}^{-1}}{\left( (\frac{\mu}{\cos\psi} - \mathrm{i}\sigma)^{2} + (\left| \boldsymbol{\xi} \right| \tan\psi)^{2} \right)^{\frac{2k+1}{2}}} \left( \left( \frac{\mu}{\cos\psi} - \mathrm{i}\sigma \right)^{2} + (\left| \boldsymbol{\xi} \right| \tan\psi)^{2} \right)^{2k} \widehat{A_{\mu}g} (\boldsymbol{\xi}, \sigma) \\
    &= \beta_{2k}^{-1} \left( \left( \frac{\mu}{\cos\psi} - \mathrm{i}\sigma \right)^{2} + (\left| \boldsymbol{\xi} \right| \tan\psi)^{2} \right)^{\frac{2k-1}{2}} \dfrac{\beta_{2k}}{\left( (\frac{\mu}{\cos\psi} - \mathrm{i}\sigma)^{2} + (\left| \boldsymbol{\xi} \right| \tan\psi)^{2} \right)^{\frac{2k-1}{2}}} \ \widehat{g} (\boldsymbol{\xi}, \sigma) \\
    &= \widehat{g} (\boldsymbol{\xi}, \sigma).
\end{align*}

\subsection{Proof of Theorem \ref{thm: odd dim of Amu}}
The proof of Theorem \ref{thm: odd dim of Amu} is more straightforward than that of Theorems \ref{thm: odd dim of Cmu} and \ref{thm: even dim of Cmu}.
For $f \in C_{c}^{\infty}(\mathbb{R}^{2k-1} \times \mathbb{R})$, suppose $g = A_{\mu}f$.
By definition of $A_{\mu}$, $g$ is smooth. Since $f$ has compact support, $g$ is bounded, and the Fourier transform of $L^{k-1}g$ is given by:
\begin{align*}
    \mathcal{F}[L^{k-1}g] (\boldsymbol{\xi}, \sigma)
    &= \left[ \left( \frac{\mu}{\cos\psi} - \mathrm{i} \sigma \right)^{2} + (\left| \boldsymbol{\xi} \right| \tan\psi)^{2} \right]^{k-1} \widehat{g} (\boldsymbol{\xi}, \sigma) \\
    &= \left[ \left( \frac{\mu}{\cos\psi} - \mathrm{i} \sigma \right)^{2} + (\left| \boldsymbol{\xi} \right| \tan\psi)^{2} \right]^{k-1} \widehat{A_{\mu}f} (\boldsymbol{\xi}, \sigma) \\
    &= \left[ \left( \frac{\mu}{\cos\psi} - \mathrm{i} \sigma \right)^{2} + (\left| \boldsymbol{\xi} \right| \tan\psi)^{2} \right]^{k-1} \dfrac{\beta_{2k-1}}{\left[ \left( \frac{\mu}{\cos\psi} - \mathrm{i} \sigma \right)^{2} + (\left| \boldsymbol{\xi} \right| \tan\psi)^{2} \right]^{k-1}} \widehat{f} (\boldsymbol{\xi}, \sigma) \\
    &= \beta_{2k-1} \widehat{f} (\boldsymbol{\xi}, \sigma).
\end{align*}
Hence, we obtain $L^{k-1}g = \beta_{2k-1}f$, which implies that $L^{k-1}g$ has compact support.

Conversely, suppose $L^{k-1}g$ has compact support.
Let's define $f$ as
\begin{equation*}
    f(\mathbf{x}, z) := \beta_{2k-1}^{-1} L^{k-1}g(\mathbf{x}, z).
\end{equation*}
Since $g$ is smooth, so is $f$.
Additionally, $f$ has compact support, because $L^{k-1}g$ has compact support.
Furthermore, we have:
\begin{align*}
    \widehat{f} (\boldsymbol{\xi}, \sigma)
    &= \beta_{2k-1}^{-1} \mathcal{F}[L^{k-1}g] (\boldsymbol{\xi}, \sigma) \\
    &= \beta_{2k-1}^{-1} \left[ \left( \frac{\mu}{\cos\psi} - \mathrm{i} \sigma \right)^{2} + (\left| \boldsymbol{\xi} \right| \tan\psi)^{2} \right]^{k-1} \widehat{g} (\boldsymbol{\xi}, \sigma).
\end{align*}
Then, the Fourier transform of $A_{\mu}f$ is as follows:
\begin{align*}
    \widehat{A_{\mu}f} (\boldsymbol{\xi}, \sigma)
    &= \dfrac{\beta_{2k-1}}{\left[ \left( \frac{\mu}{\cos\psi} - \mathrm{i} \sigma \right)^{2} + (\left| \boldsymbol{\xi} \right| \tan\psi)^{2} \right]^{k-1}} \widehat{f} (\boldsymbol{\xi}, \sigma) \\
    &= \dfrac{\beta_{2k-1}}{\left[ \left( \frac{\mu}{\cos\psi} - \mathrm{i} \sigma \right)^{2} + (\left| \boldsymbol{\xi} \right| \tan\psi)^{2} \right]^{k-1}} \beta_{2k-1}^{-1} \left[ \left( \frac{\mu}{\cos\psi} - \mathrm{i} \sigma \right)^{2} + (\left| \boldsymbol{\xi} \right| \tan\psi)^{2} \right]^{k-1} \widehat{g} (\boldsymbol{\xi}, \sigma) \\
    &= \widehat{g} (\boldsymbol{\xi}, \sigma).
\end{align*}

\subsection{Proof of Theorem \ref{thm: even dim of Amu}}
Similar to the proof of Theorem \ref{thm: odd dim of Amu}, for $f \in C_{c}^{\infty}(\mathbb{R}^{2k} \times \mathbb{R})$ suppose $g = A_{\mu}f$. By definition of $A_{\mu}$, $g$ is smooth. Since $f$ has compact support, $g$ is bounded and
\begin{align*}
    \mathcal{F}[L^{2k-1} A_{\mu}g] (\boldsymbol{\xi}, \sigma)
    &= \left[ \left( \frac{\mu}{\cos\psi} - \mathrm{i} \sigma \right)^{2} + (\left| \boldsymbol{\xi} \right| \tan\psi)^{2} \right]^{2k-1} \widehat{A_{\mu}g} (\boldsymbol{\xi}, \sigma) \\
    &= \left[ \left( \frac{\mu}{\cos\psi} - \mathrm{i} \sigma \right)^{2} + (\left| \boldsymbol{\xi} \right| \tan\psi)^{2} \right]^{2k-1} \dfrac{\beta_{2k}}{\left[ \left( \frac{\mu}{\cos\psi} - \mathrm{i} \sigma \right)^{2} + (\left| \boldsymbol{\xi} \right| \tan\psi)^{2} \right]^{\frac{2k-1}{2}}} \widehat{A_{\mu}f} (\boldsymbol{\xi}, \sigma) \\
    &= \left[ \left( \frac{\mu}{\cos\psi} - \mathrm{i} \sigma \right)^{2} + (\left| \boldsymbol{\xi} \right| \tan\psi)^{2} \right]^{2k-1} \dfrac{(\beta_{2k})^{2}}{\left[ \left( \frac{\mu}{\cos\psi} - \mathrm{i} \sigma \right)^{2} + (\left| \boldsymbol{\xi} \right| \tan\psi)^{2} \right]^{2k-1}} \widehat{f} (\boldsymbol{\xi}, \sigma) \\
    &= (\beta_{2k})^{2} \widehat{f} (\boldsymbol{\xi}, \sigma).
\end{align*}
Hence, we obtain $L^{2k-1}A_{\mu}g = (\beta_{2k})^{2}f$, which implies that $L^{2k-1}A_{\mu}g$ has compact support.

Conversely, suppose $L^{2k-1}A_{\mu}g$ has compact support;
$f$ can be defined as:
\begin{equation*}
    f(\mathbf{x}, z) := \beta_{2k+1}^{-2} L^{2k-1}A_{\mu}g(\mathbf{x}, z).
\end{equation*}
Since $g$ is smooth, so is $f$.
In addition, $f$ has compact support because $L^{k}g$ has compact support.
Furthermore, we have:
\begin{align*}
    \widehat{f} (\boldsymbol{\xi}, \sigma)
    &= \beta_{2k}^{-2} \mathcal{F}[L^{2k-1}A_{\mu}g] (\boldsymbol{\xi}, \sigma) \\
    &= \beta_{2k}^{-1} \left[ \left( \frac{\mu}{\cos\psi} - \mathrm{i} \sigma \right)^{2} + (\left| \boldsymbol{\xi} \right| \tan\psi)^{2} \right]^{\frac{2k-1}{2}} \widehat{g} (\boldsymbol{\xi}, \sigma).
\end{align*}
Then the Fourier transform of $A_{\mu}f$ is as follows:
\begin{align*}
    \widehat{A_{\mu}f} (\boldsymbol{\xi}, \sigma)
    &= \dfrac{\beta_{2k}}{\left[ \left( \frac{\mu}{\cos\psi} - \mathrm{i} \sigma \right)^{2} + (\left| \boldsymbol{\xi} \right| \tan\psi)^{2} \right]^{\frac{2k-1}{2}}} \widehat{f} (\boldsymbol{\xi}, \sigma) \\
    &= \dfrac{\beta_{2k}}{\left[ \left( \frac{\mu}{\cos\psi} - \mathrm{i} \sigma \right)^{2} + (\left| \boldsymbol{\xi} \right| \tan\psi)^{2} \right]^{\frac{2k-1}{2}}} \beta_{2k}^{-1} \left[ \left( \frac{\mu}{\cos\psi} - \mathrm{i} \sigma \right)^{2} + (\left| \boldsymbol{\xi} \right| \tan\psi)^{2} \right]^{\frac{2k-1}{2}} \widehat{g} (\boldsymbol{\xi}, \sigma) \\
    &= \widehat{g} (\boldsymbol{\xi}, \sigma).
\end{align*}

\section{Conclusion}
In this paper, we studied the range conditions for the operators $C_{\mu}$ and $A_{\mu}$.
We demonstrated that differential operator $L$ is vital in describing the range conditions for $C_{\mu}$ and $A_{\mu}$.
Additionally, we showed that the range conditions for $C_{\mu}$ and $A_{\mu}$ vary depending on the dimension of the domain.
Our results contribute significantly to the understanding of the attenuated conical Radon transforms.
Furthermore, in a future study, we plan to address the inversion and range description of generalized attenuated conical Radon transforms.

\section*{Acknowledgment}
This work was supported by the National Research Foundation of Korea (NRF-2022R1C1C1003464 and RS-2023-00217116).
Professor Sunghwan Moon provided valuable ideas and comments.

\section*{Appendix}
\renewcommand{\thesubsection}{\Alph{subsection}}

In this Appendix, we provide the proof of Lemma \ref{lem: fourier transform of Cmuf and Amuf}. The following identity is useful for the proof \cite[p198]{natterer2001mathematics}:
\begin{equation}\label{eq: Funk-Hecke}
    \int\limits_{S^{n-1}}e^{-\mathrm{i}\sigma \theta \cdot \boldsymbol{\xi}} \mathrm{d}\boldsymbol{\xi} = (2\pi)^{\frac{n}{2}}\sigma^{\frac{2-n}{2}} J_{\frac{n-2}{2}}(\sigma).
\end{equation}

\subsection*{Proof of Lemma \ref{lem: fourier transform of Cmuf and Amuf}}
According to \cite{gouia2018inversion}, we have
\begin{equation*}
    \widehat{C_{\mu}f} (\boldsymbol{\xi}, \sigma) = \dfrac{(2\pi)^{\frac{n}{2}} \left| \boldsymbol{\xi} \right|^{1-\frac{n}{2}} \tan^{\frac{n}{2}}\psi}{\cos\psi} \int\limits_{0}^{\infty} J_{\frac{n-2}{2}}(z \left| \boldsymbol{\xi} \right| \tan\psi) z^{\frac{n}{2}}e^{-(\frac{\mu}{\cos\psi} - \mathrm{i}\sigma)z} \mathrm{d}z \widehat{f} (\boldsymbol{\xi}, \sigma),
\end{equation*}
where $J_{\nu}$ is the Bessel function of the first kind of order $\nu$.
The integral of the above equation is given by
\begin{align*}
    \int\limits_{0}^{\infty} J_{\frac{n-2}{2}}(z \left| \boldsymbol{\xi} \right| \tan\psi) z^{\frac{n}{2}}e^{-z(\frac{\mu}{\cos\psi} - \mathrm{i}\sigma)} \mathrm{d}z
    &= (\left| \boldsymbol{\xi} \right| \tan\psi)^{-\frac{1}{2}}\int\limits_{0}^{\infty} z^{\frac{n-1}{2}}e^{-(\frac{\mu}{\cos\psi} - \mathrm{i}\sigma)z} J_{\frac{n-2}{2}}(z \left| \boldsymbol{\xi} \right| \tan\psi) (z \left| \boldsymbol{\xi} \right| \tan\psi)^{\frac{1}{2}} \mathrm{d}z \\ 
    &= (\left| \boldsymbol{\xi} \right| \tan\psi)^{-\frac{1}{2}} \dfrac{\pi^{-\frac{1}{2}}2^{\frac{n}{2}} \Gamma(\frac{n+1}{2}) (\frac{\mu}{\cos\psi} - \mathrm{i}\sigma) (\left| \boldsymbol{\xi} \right| \tan\psi)^{\frac{n-1}{2}}}{\left((\frac{\mu}{\cos\psi} - \mathrm{i}\sigma)^{2} + (\left| \boldsymbol{\xi} \right| \tan\psi)^{2} \right)^{\frac{n+1}{2}}} \\ 
    &= \dfrac{\pi^{-\frac{1}{2}}2^{\frac{n}{2}} \Gamma(\frac{n+1}{2}) (\frac{\mu}{\cos\psi} - \mathrm{i}\sigma) (\left| \boldsymbol{\xi} \right| \tan\psi)^{\frac{n-2}{2}}}{\left((\frac{\mu}{\cos\psi} - \mathrm{i}\sigma)^{2} + (\left| \boldsymbol{\xi} \right| \tan\psi)^{2} \right)^{\frac{n+1}{2}}}.
\end{align*}
where in the second line, we use the following identity \cite[p198]{natterer2001mathematics}: For $\operatorname{Re}\ a >0$ and $\operatorname{Re}\ \nu >-1$,
\begin{equation}
    \int\limits_{0}^{\infty} x^{\nu+\frac{1}{2}}e^{-ax} J_{\nu}(xy)(xy)^{\frac{1}{2}}dx
    = \dfrac{\pi^{-\frac{1}{2}}2^{\nu +1 } \Gamma(\nu + \frac{3}{2}) a y^{\nu + \frac{1}{2}}}{(a^{2} + y^{2})^{\nu + \frac{3}{2}}},\quad y > 0.
\end{equation}

Therefore, we obtain
\begin{align*}
    \widehat{C_{\mu}f} (\boldsymbol{\xi}, \sigma)
    &= \dfrac{(2\pi)^{\frac{n}{2}} \left| \boldsymbol{\xi} \right|^{1-\frac{n}{2}} \tan^{\frac{n}{2}}\psi}{\cos\psi} \dfrac{\pi^{-\frac{1}{2}}2^{\frac{n}{2}} \Gamma(\frac{n+1}{2}) (\frac{\mu}{\cos\psi} - \mathrm{i}\sigma) (\left| \boldsymbol{\xi} \right| \tan\psi)^{\frac{n-2}{2}}}{\left((\frac{\mu}{\cos\psi} - \mathrm{i}\sigma)^{2} + (\left| \boldsymbol{\xi} \right| \tan\psi)^{2} \right)^{\frac{n+1}{2}}} \widehat{f} (\boldsymbol{\xi}, \sigma) \\ 
    &= \dfrac{2^{n}\pi^{\frac{n-1}{2}} \Gamma(\frac{n+1}{2}) \tan^{n-1}\psi }{\cos\psi} \dfrac{\frac{\mu}{\cos\psi} - \mathrm{i}\sigma}{\left((\frac{\mu}{\cos\psi} - \mathrm{i}\sigma)^{2} + (\left| \boldsymbol{\xi} \right| \tan\psi)^{2} \right)^{\frac{n+1}{2}}} \widehat{f} (\boldsymbol{\xi}, \sigma) \\ 
    &= \alpha_{n} \dfrac{\frac{\mu}{\cos\psi} - \mathrm{i}\sigma}{\left((\frac{\mu}{\cos\psi} - \mathrm{i}\sigma)^{2} + (\left| \boldsymbol{\xi} \right| \tan\psi)^{2} \right)^{\frac{n+1}{2}}} \widehat{f} (\boldsymbol{\xi}, \sigma).
\end{align*}
The proof for $A_{\mu}f$ is similar. The Fourier transform of $A_{\mu}f$ is given by 
\begin{align*}
    \widehat{A_{\mu}} (\boldsymbol{\xi}, \sigma)
    &= \dfrac{\tan^{n-1}\psi}{\cos\psi} \int\limits_{\mathbb{R}} \int\limits_{\mathbb{R}^{n}} \int\limits_{S^{n-1}} \int\limits_{0}^{\infty} f\left( (\mathbf{u}, v) + (z \tan \psi \boldsymbol{\omega}, z) \right)e^{-\mathrm{i} \mathbf{u} \cdot \boldsymbol{\xi}} e^{-\mathrm{i}v \sigma} e^{-\frac{\mu}{\cos\psi}z} z^{n-2} \mathrm{d}z \mathrm{d}S(\boldsymbol{\omega}) \mathrm{d}\mathbf{u} \mathrm{d}v \\\ 
    &= \dfrac{\tan^{n-1}\psi}{\cos\psi} \int\limits_{S^{n-1}} \int\limits_{0}^{\infty} e^{-\mathrm{i} z \tan \psi \boldsymbol{\omega} \cdot \boldsymbol{\xi}} e^{-\mathrm{i} \sigma z} e^{-\frac{\mu}{\cos\psi}z} z^{n-2} \mathrm{d}z \mathrm{d}S(\boldsymbol{\omega}) \widehat{f}(\boldsymbol{\xi}, \sigma) \\\ 
    &= \dfrac{(2\pi)^{\frac{n}{2}}\tan^{\frac{n}{2}}\psi \left| \boldsymbol{\xi} \right|^{1-\frac{n}{2}}}{\cos\psi} \int\limits_{0}^{\infty} J_{\frac{n-2}{2}}(z \left| \boldsymbol{\xi} \right| \tan\psi) e^{-(\frac{\mu}{\cos\psi} - \mathrm{i} \sigma)z} z^{\frac{n-2}{2}} \mathrm{d}z \widehat{f} (\boldsymbol{\xi}, \sigma),
\end{align*}
where in the last line, we use identity \eqref{eq: Funk-Hecke}.
Moreover, we get
\begin{align*}
    \int\limits_{0}^{\infty} J_{\frac{n-2}{2}}(z \left| \boldsymbol{\xi} \right| \tan\psi) e^{-(\frac{\mu}{\cos\psi} - \mathrm{i} \sigma)z} z^{\frac{n-2}{2}} \mathrm{d}z
    &= (\left| \boldsymbol{\xi} \right| \tan\psi)^{-\frac{1}{2}} \int\limits_{0}^{\infty} z^{\frac{n-3}{2}} e^{-(\frac{\mu}{\cos\psi} - \mathrm{i} \sigma)z} J_{\frac{n-2}{2}}(z \left| \boldsymbol{\xi} \right| \tan\psi) (z \left| \boldsymbol{\xi} \right| \tan\psi)^{\frac{1}{2}} \mathrm{d}z \\\ 
    &= (\left| \boldsymbol{\xi} \right| \tan\psi)^{-\frac{1}{2}} \dfrac{\pi^{-\frac{1}{2}}2^{\frac{n-2}{2}} \Gamma(\frac{n-1}{2}) (\left| \boldsymbol{\xi} \right| \tan\psi)^{\frac{n-1}{2}}}{\left( (\frac{\mu}{\cos\psi} - \mathrm{i} \sigma)^{2} + (\left| \boldsymbol{\xi} \right| \tan\psi)^{2} \right)^{\frac{n-1}{2}}} \\\ 
    &= \dfrac{\pi^{-\frac{1}{2}}2^{\frac{n-2}{2}} \Gamma(\frac{n-1}{2}) (\left| \boldsymbol{\xi} \right| \tan\psi)^{\frac{n-2}{2}}}{\left( (\frac{\mu}{\cos\psi} - \mathrm{i} \sigma)^{2} + (\left| \boldsymbol{\xi} \right| \tan\psi)^{2} \right)^{\frac{n-1}{2}}}.
\end{align*}
where in the second line, we use the following identity \cite[p29]{bateman1954tables}: For $\operatorname{Re} a > 0$ and $\operatorname{Re} \nu > -\frac{1}{2}$ ,
\begin{equation*}
    \int\limits_{0}^{\infty} x^{\nu - \frac{1}{2}}e^{-ax} J_{\nu}(xy)(xy)^{\frac{1}{2}}dx
    = \dfrac{\pi^{-\frac{1}{2}}2^{\nu} \Gamma(\nu + \frac{1}{2}) y^{\nu + \frac{1}{2}}}{(a^{2} + y^{2})^{\nu + \frac{1}{2}}},\quad y > 0.
\end{equation*}
Therefore, the Fourier transform of $A_{\mu}f$ is given by
\begin{align*}
    \widehat{A_{\mu}f} (\boldsymbol{\xi}, \sigma)
    &= \dfrac{(2\pi)^{\frac{n}{2}}\tan^{\frac{n}{2}}\psi \left| \boldsymbol{\xi} \right|^{1-\frac{n}{2}}}{\cos\psi} \dfrac{\pi^{-\frac{1}{2}}2^{\frac{n-2}{2}} \Gamma(\frac{n-1}{2}) (\left| \boldsymbol{\xi} \right| \tan\psi)^{\frac{n-2}{2}}}{\left( (\frac{\mu}{\cos\psi} - \mathrm{i} \sigma)^{2} + (\left| \boldsymbol{\xi} \right| \tan\psi)^{2} \right)^{\frac{n-1}{2}}}  \widehat{f} (\boldsymbol{\xi}, \sigma) \\\ 
    &= \dfrac{2^{n-1}\pi^{\frac{n-1}{2}}\Gamma(\frac{n-1}{2})\tan^{n-1}\psi}{\cos\psi\left( (\frac{\mu}{\cos\psi} - \mathrm{i} \sigma)^{2} + (\left| \boldsymbol{\xi} \right| \tan\psi)^{2} \right)^{\frac{n-1}{2}}}  \widehat{f} (\boldsymbol{\xi}, \sigma) \\\ 
    &= \dfrac{\beta_{n}}{\left( (\frac{\mu}{\cos\psi} - \mathrm{i} \sigma)^{2} + (\left| \boldsymbol{\xi} \right| \tan\psi)^{2} \right)^{\frac{n-1}{2}}}  \widehat{f} (\boldsymbol{\xi}, \sigma).
\end{align*}

\bibliographystyle{unsrt}
\bibliography{ref}

\end{document}